\documentclass[a4paper,11pt]{article}
\usepackage{amsfonts}
\usepackage{amstext}
\usepackage{amsthm}
\usepackage{amsmath}
\usepackage{amssymb}
\usepackage{latexsym}
\usepackage{graphicx}
\usepackage[latin1]{inputenc}
\newcommand{\R}{\Bbb{R}}
\newcommand{\N}{\Bbb{N}}

\newtheorem{teor}{Theorem}[section]
\newtheorem{propo}{Proposition}[section]

\newtheorem{rem}{Remark}[section]

\newcommand{\n}{\noindent}

\newcommand {\fim}{\rule{0.5em}{0.5em}}

\begin{document}

\title{On positive viscosity solutions of fractional\\ Lane-Emden systems
}

\author{\textbf{Edir Junior Ferreira Leite \footnote{\textit{E-mail addresses}:
edirjrleite@ufv.br (E.J.F. Leite)}}\\ {\small\it Departamento de Matem\'{a}tica,
Universidade Federal de Vi\c{c}osa,}\\ {\small\it CCE, 36570-000, Vi\c{c}osa, MG, Brazil}\\
\textbf{Marcos Montenegro \footnote{\textit{E-mail addresses}:
montene@mat.ufmg.br (M.
Montenegro)}}\\ {\small\it Departamento de Matem\'{a}tica,
Universidade Federal de Minas Gerais,}\\ {\small\it Caixa Postal
702, 30123-970, Belo Horizonte, MG, Brazil}}

\date{}{

\maketitle

\markboth{abstract}{abstract}
\addcontentsline{toc}{chapter}{abstract}

\hrule \vspace{0,2cm}

\n {\bf Abstract}

In this paper we discuss on existence, nonexistence and uniqueness of positive viscosity solution for the following coupled system involving fractional Laplace operators on a smooth bounded domain $\Omega$ in $\R^n$:

\[
\left\{
\begin{array}{llll}
(-\Delta)^{s}u = v^p & {\rm in} \ \ \Omega\\
(-\Delta)^{s}v = u^q & {\rm in} \ \ \Omega\\
u= v=0 & {\rm in} \ \ \R^n\setminus\Omega
\end{array}
\right.
\]

\n By mean of an appropriate variational framework and a H\"{o}lder type regularity result for suitable weak solutions of the above system, we prove that such a system admits at least one positive viscosity solution for any $0 < s < 1$, provided that $p, q > 0$, $pq \neq 1$ and the couple $(p,q)$ is below the critical hyperbole

\[
\frac{1}{p + 1} + \frac{1}{q + 1} = \frac{n - 2s}{n}
\]
whenever $n > 2s$.

\n Moreover, by using maximum principles for the fractional Laplace operator, we show that uniqueness occurs whenever $pq < 1$. Lastly, assuming $\Omega$ is star-shaped, by using a Rellich type variational identity, we prove that no such a solution exists if $(p,q)$ is on or above the critical hyperbole. A crucial point in our proofs is proving, given a critical point $u \in W_0^{s, \frac{p+1}{p}}(\Omega) \cap W^{2s, \frac{p+1}{p}}(\Omega)$ of a related functional, that there is a function $v$ in an appropriate Sobolev space (Proposition 2.1) so that $(u,v)$ is a weak solution of the above system and a bootstrap argument can be applied successfully in order to establish its H\"{o}lder regularity (Proposition 3.1). The difficulty is caused mainly by the absence of a general Agmon-Zygmund theory for $0 < s <1$. In particular, the employed idea is new and different of the corresponding one for $s = 1$.

\vspace{0.5cm}
\hrule\vspace{0.2cm}

\section{Introduction and main results}

This work is devoted to the study of existence, uniqueness and nonexistence of positive viscosity solutions for nonlocal elliptic systems on bounded domains which will be described henceforth.

The fractional Laplace operator (or fractional Laplacian) of order $2s$, with $0 < s < 1$, denoted by $(-\Delta)^{s}$, is defined as

\[
(-\Delta)^{s}u(x) = C(n,s)\, {\rm P.V.}\int\limits_{\R^{n}}\frac{u(x)-u(y)}{\vert x-y\vert^{n+2s}}\; dy\, ,
\]
or equivalently,

\[
(-\Delta)^{s}u(x)=-\frac{1}{2}C(n,s)\int\limits_{\R^{n}}\frac{u(x+y)+u(x-y)-2u(x)}{|y|^{n+2s}}\; dy
\]

\n for all $x \in \R^{n}$, where P.V. denotes the principal value of the first integral and

\[
C(n,s) = \left(\int\limits_{\R^{n}}\frac{1-\cos(\zeta_{1})}{\vert\zeta\vert^{n+2s}}\; d\zeta\right)^{-1}
\]

\n with $\zeta = (\zeta_1, \ldots, \zeta_n) \in \R^n$.

Remark that $(-\Delta)^{s}$ is a nonlocal operator on functions compactly supported in $\R^n$. Moreover, for any function $u \in C^{\infty}_{0}(\R^{n})$,

\[
\lim_{s\rightarrow 1^{-}}(-\Delta)^{s} u(x) = -\Delta u(x)
\]

\n for all $x \in \R^n$, so that the operator $(-\Delta)^{s}$ interpolates the Laplace operator in $\R^n$. Factional Laplace operators arise naturally in different areas such as Probability, Finance, Physics, Chemistry and Ecology, see \cite{A}.

A closely related operator, but different from $(-\Delta)^{s}$, is the spectral fractional Laplace operator $\mathcal{A}^{s}$ which is defined in terms of the Dirichlet spectra of the Laplace operator on $\Omega$. Roughly, for a $L^2$-orthonormal basis of eigenfunctions $(\varphi_k)$ corresponding to eigenvalues $(\lambda_k)$ of the Laplace operator with zero Dirichlet boundary values on $\partial \Omega$, the operator $\mathcal{A}^s$ is defined as $\mathcal{A}^{s} u = \sum_{k=1}^\infty c_k \lambda_k^s \varphi_k$, where $c_k$, $k \geq 1$, are the coefficients of the expansion $u = \sum_{k=1}^\infty c_k \varphi_k$.

After the work \cite{CaSi} on the characterization for any $0 < s < 1$ of the operator $(-\Delta)^{s}$ in terms of a Dirichlet-to-Neumann map associated to a suitable extension problem, a great deal of attention has been dedicated in the last years to nonlinear nonlocal problems of the kind

\begin{equation}\label{3}
\left\{
\begin{array}{rrll}
(-\Delta)^{s} u &=& f(x,u) & {\rm in} \ \ \Omega\\
u &=& 0  & {\rm in} \ \ \R^n\setminus\Omega
\end{array}
\right.
\end{equation}
where $\Omega$ is a smooth bounded open subset of $\R^{n}$, $n \geq 1$ and $0 < s < 1$.

Several works have been focused on the existence \cite{ros131, ros137, ros138, auto2, auto1, ros215, niang, ros264, ros265, ros267, val}, nonexistence \cite{ros129, val}, symmetry \cite{ros21, ros95} and regularity \cite{ros10, cabre1, ROS} of viscosity solutions, among other qualitative properties \cite{ros1, ros130}. For developments related to (\ref{3}) involving the spectral fractional Laplace operator $\mathcal{A}^{s}$, we refer to \cite{BaCPS, BrCPS, CT, CDDS, choi, CK, tan1, T, yang}, among others.

A specially important example is given by the power function $f(x,u) = u^p$ for $p > 0$, in which case (\ref{3}) is called fractional Lane-Emden problem. Recently, it has been proved in \cite{SV} that this problem admits at least one positive viscosity solution for $1 < p < \frac{n + 2s}{n - 2s}$. The nonexistence has been established in \cite{ROS2} whenever $p \geq \frac{n + 2s}{n - 2s}$ and $\Omega$ is star-shaped. These results were known long before for $s = 1$, see the classical references \cite{AR, GS, ros237, Rabinowitz} and the survey \cite{Pucci}.

We here are interested in studying the following vector counterpart of the fractional Lane-Emden problem:

\begin{equation} \label{1}
\left\{
\begin{array}{llll}
(-\Delta)^{s}u = v^p & {\rm in} \ \ \Omega\\
(-\Delta)^{s}v = u^q & {\rm in} \ \ \Omega\\
u= v=0 & {\rm in} \ \ \R^n\setminus\Omega
\end{array}
\right.
\end{equation}

\n for $p, q > 0$, which inspires the title of this work.

For $s = 1$, the problem (\ref{1}) and a number of its generalizations have been widely investigated in the literature, see for instance the survey \cite{DG} and references therein. Specifically, notions of sublinearity, superlinearity and criticality (subcriticality, supercriticality) have been introduced in \cite{FM, Mi1, Mi2, SZ1}. In fact, the behavior of (\ref{1}) is sublinear when $pq < 1$, superlinear when $pq > 1$ and critical (subcritical, supercritical) when $n \geq 3$ and $(p,q)$ is on (below, above) the hyperbole, known as critical hyperbole,

\[
\frac{1}{p+1}+\frac{1}{q+1}=\frac{n-2}{n}\, .
\]

\n When $pq = 1$, its behavior is resonant and the corresponding eigenvalue problem has been addressed in \cite{marcos}. The sublinear case has been studied in \cite{FM} where the existence and uniqueness of positive classical solution is proved. The superlinear-subcritical case has been covered in the works \cite{CFM}, \cite{DF}, \cite{DR} and \cite{vander} where the existence of at least one positive classical solution is derived. Lastly, the nonexistence of positive classical solutions has been established in \cite{Mi1} on star-shaped domains.

In this work we discuss the existence and nonexistence of positive viscosity solution of (\ref{1}) for any $0 < s < 1$. We determine the precise set of exponents $p$ and $q$ for which the problem (\ref{1}) admits always a positive viscosity solution. In particular, we extend the above-mentioned results corresponding to the fractional Lane-Emden problem for $0 < s \leq 1$ and to the Lane-Emden system involving the Laplace operator. As a byproduct, the notions of sublinearity, superlinearity and criticality (subcriticality, supercriticality) to the problem (\ref{1}) are naturally extended for any $0 < s < 1$.

The ideas involved in our proofs base on variational methods, $C^\beta$ regularity of weak solutions and an integral variational identity satisfied by positive viscosity solutions of (\ref{1}). We shall introduce a variational framework in order to establish the existence of nontrivial nonnegative weak solution of (\ref{1}) in a suitable sense. In our formulation, the function $u$ arises as a nonzero critical point of a functional defined on an appropriate space of functions and then, in a natural way, we construct a function $v$ so that the couple $(u,v)$ is a weak solution of (\ref{1}). The construction of $v$ is a strategic step in our approach (see Section 2), once we do not have a Calderón-Zygmund theory available for fractional operators. Using the $C^\beta$ regularity result (to be proved in Section 3) for weak solutions of (\ref{1}) and maximum principles for fractional Laplace operators, we then deduce that the couple $(u,v)$ is a positive viscosity solution of (\ref{1}). Moreover, we prove its uniqueness when $pq < 1$. The key tool used in the nonexistence proof is a Rellich type variational identity (to be proved in Section 4) to positive viscosity solutions of (\ref{1}). The proof of the $C^\beta$ regularity consists in first showing that weak solutions of (\ref{1}) belong to $L^r(\Omega) \times L^r(\Omega)$ for every $r \geq 1$ and then applying to each equation the $C^\beta$ regularity up to the boundary proved recently in \cite{ROS1}. The proof of the variational identity to Lane-Emden systems uses the Pohozaev variational identity to fractional elliptic equations obtained recently in \cite{ROS2}.

Other questions have recently been discussed in some papers which mention ours. We quote for example the works \cite{ChK} on asymptotic behavior of minimal energy solutions, \cite{YP} on symmetry of solutions and \cite{ZY} on Liouville type theorems on half spaces for fractional systems. Related systems also have been investigated by using other methods. We refer to the work \cite{EM} for systems involving different operators $(-\Delta)^{s}$ and $(-\Delta)^{t}$ in each one of equations and to the work \cite{choi} for systems involving the spectral fractional operator $\mathcal{A}^{s}$.

In order to state our three main theorems, we should first introduce the concept of positive viscosity solution to (\ref{1}). A couple $(u,v)$ of continuous functions in $\R^n$ is said to be a viscosity subsolution (supersolution) of (\ref{1}) if each point $x_{0} \in \Omega$ admits a neighborhood $U$ with $\overline{U} \subset \Omega$ such that for any $\varphi,\psi \in C^{2}(\overline{U})$ satisfying $u(x_{0}) = \varphi(x_{0})$, $v(x_{0}) = \psi(x_{0})$, $u \geq \varphi$ and $v \geq \psi$ in $U$, the functions

\begin{equation} \label{2.2}
\overline{u}=\left\{\begin{array}{lllc}
\varphi & {\rm in} \ U  \\
u & {\rm in} \ \R^{n}\setminus U
\end{array}\right. \ \ {\rm and}\ \
\overline{v}=\left\{\begin{array}{lllc}
\psi & {\rm in} \ U  \\
v & {\rm in} \ \R^{n}\setminus U
\end{array}\right.
\end{equation}

\n satisfy
\[
(-\Delta)^{s}\overline{u}(x_{0}) \leq\ (\geq)\ |v(x_{0})|^{p - 1} v(x_{0}) \ \ {\rm and}\ \
(-\Delta)^{s}\overline{v}(x_{0}) \leq\ (\geq)\ |u(x_{0})|^{q - 1} u(x_{0})\, .
\]

\n A couple $(u,v)$ of functions is said to be a viscosity solution of (\ref{1}) if it is simultaneously a viscosity subsolution and supersolution. If further $u$ and $v$ are positive in $\Omega$ and nonnegative in $\R^n$, we say that $(u,v)$ is a positive viscosity solution.

\begin{teor}\label{teo1} (sublinear case)
Let $\Omega$ be a smooth bounded open subset of $\R^n$, $n \geq 1$ and $0 < s < 1$. Assume that $p,q>0$ and $pq < 1$. Then the problem (\ref{1}) admits a unique positive viscosity solution.
\end{teor}

\begin{teor}\label{teo2} (superlinear-subcritical case)
Let $\Omega$ be a smooth bounded open subset of $\R^n$, $n \geq 1$ and $0 < s < 1$. Assume that $p,q>0$, $pq > 1$ and

\begin{equation}\label{5}
\frac{1}{p + 1} + \frac{1}{q + 1} > \frac{n - 2s}{n}\, .
\end{equation}

\n whenever $n > 2s$. Then the problem (\ref{1}) admits at least one positive viscosity solution.
\end{teor}

\begin{teor}\label{teo3} (critical and supercritical cases)
Let $\Omega$ be a smooth bounded open subset of $\R^n$, $n > 2s$ and $0 < s < 1$. Assume that $\Omega$ is star-shaped, $p,q>0$ and

\begin{equation}\label{4}
\frac{1}{p + 1} + \frac{1}{q + 1} \leq \frac{n - 2s}{n}\, .
\end{equation}

\n Then the problem (\ref{1}) admits no positive viscosity solution.
\end{teor}

\begin{figure}[ht]
\centering
\includegraphics[scale=1.2]{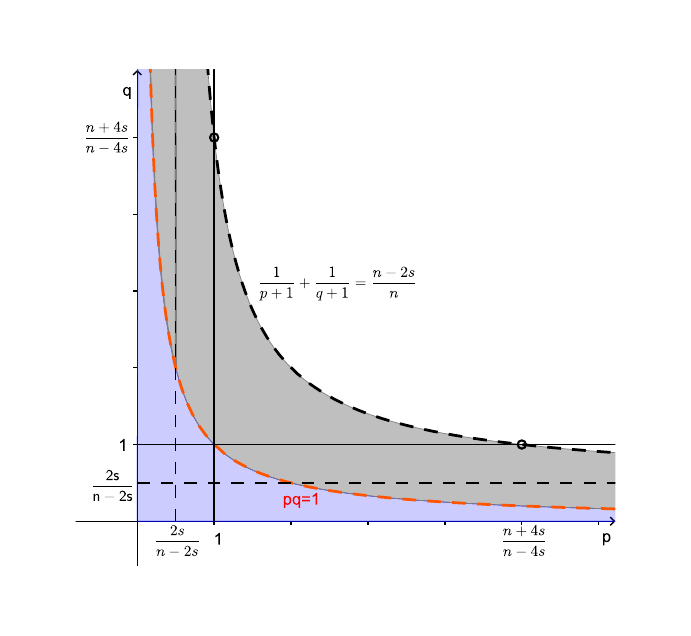}
\caption{The existence range of couples $(p,q)$ when $n>4s$.}
\end{figure}

For dimensions $n > 2s$, the hyperbole

\begin{equation}\label{6}
\frac{1}{p + 1} + \frac{1}{q + 1} = \frac{n - 2s}{n}
\end{equation}
is called critical hyperbole associated to the Lane-Emden system (\ref{1}). It appears naturally in the context of Sobolev embedding theorems and integral variational identities. Note also that the curve $(p,q)$ given by the hyperbole $pq = 1$ splits the behavior of (\ref{1}) into sublinear and superlinear one.

The remainder of paper is organized into six sections. In Section 2, we briefly recall some definitions and facts dealing with fractional Sobolev spaces and introduce the variational framework and the adequate concept of weak solution to be used in the existence proofs. In Section 3, we prove the $C^\beta$ regularity of weak solutions of (\ref{1}) into the subcritical context according to the hyperbole (\ref{6}). In Section 4, we establish a Rellich variational identity to positive viscosity solutions of (\ref{1}). In Section 5, we prove Theorem \ref{teo1} by using a direct minimization approach, the regularity provided in the third section and maximum principles. In section 6, we prove Theorem \ref{teo2} by using the mountain pass theorem and the same regularization result. Finally, in Section 7, we prove Theorem \ref{teo3} by applying the variational identity obtained in Section 4.\\

\section{Preliminaries and variational setting}

In this section, we recall the definition of fractional Sobolev spaces on bounded open subsets of $\R^n$ and present the variational formulation to be used in the proofs of Theorems \ref{teo1} and \ref{teo2}.

We start by fixing a parameter $0 < s < 1$. Let $\Omega$ be an open subset of $\R^n$ with $n \geq 1$. For any $r \in (1,+\infty)$, one defines the fractional Sobolev space $W^{s,r}(\Omega)$ as

\begin{equation}\label{frac2.1}
W^{s,r}(\Omega):=\left\{u\in L^r(\Omega): \frac{\vert u(x)-u(y)\vert}{\vert x-y\vert^{\frac{n}{r}+s}}\in L^r(\Omega\times\Omega)\right\}\, ,
\end{equation}
that is, an intermediary Banach space between $L^r(\Omega)$ and $W^{1,r}(\Omega)$ induced with the norm
\begin{equation}\label{frac2.2}
\Vert u\Vert_{W^{s,r}(\Omega)}:=\left(\int\limits_{\Omega}\vert u\vert^r dx + \int\limits_{\Omega}\int\limits_{\Omega}\frac{\vert u(x)-u(y)\vert^r}{\vert x-y\vert^{n+sr}}dxdy\right)^{\frac{1}{r}}\, ,
\end{equation}
where the term
\[
[u]_{W^{s,r}(\Omega)}:=\left(\int\limits_{\Omega}\int\limits_{\Omega}\frac{\vert u(x)-u(y)\vert^r}{\vert x-y\vert^{n+sr}}dxdy\right)^{\frac{1}{r}}
\]
is called Gagliardo semi-norm of $u$.

Let $s \in \R \setminus \N$ with $s\geq 1$. The space $W^{s,r}(\Omega)$ is defined as
\[
W^{s,r}(\Omega)=\{u\in W^{[s],r}(\Omega) : D^ju\in W^{s-[s],r}(\Omega),\forall j, \vert j\vert=[s]\}\, ,
\]

\n where $[s]$ is the largest integer smaller than $s$, $j$ denotes the $n$-uple $(j_1, \ldots, j_n) \in \N^n$ and $|j|$ denotes the sum $j_1 + \ldots + j_n$.

It is clear that $W^{s,r}(\Omega)$ endowed with the norm

\begin{equation}\label{frac2.11}
\Vert u\Vert_{W^{s,r}(\Omega)}=\left(\Vert u\Vert^r_{W^{[s],r}(\Omega)}+ [u]^r_{W^{s-[s],r}(\Omega)}\right)^{\frac{1}{r}}
\end{equation}
is a reflexive Banach space.

Clearly, if $s=m$ is an integer, the space $W^{s,r}(\Omega)$ coincides with the Sobolev space $W^{m,r}(\Omega)$.

Let $W^{s,r}_0(\Omega)$ denote the closure of $C^\infty_0(\Omega)$ with respect to the norm $\Vert\cdot\Vert_{W^{s,r}(\Omega)}$ defined in (\ref{frac2.11}). For $0 < s \leq 1$, we have
\[
W^{s,r}_0(\Omega)=\{u\in W^{s,r}(\R^n) : u=0 \text{ in }\R^n\setminus\Omega\}\,
\]
and $W^{s,2}_0(\Omega)=H^s_0(\Omega)$. For more details on the above claims, we refer to \cite{tartar}.

We are ready to introduce the variational framework associated to (\ref{1}).

Let $\Omega$ be a smooth bounded open subset of $\R^n$, $n \geq 1$ and $0 < s < 1$. In order to inspire our formulation, assume that the couple $(u,v)$ of nonnegative functions is roughly a solution of (\ref{1}). From the first equation, we have $v = \left((-\Delta)^{s} u\right)^{\frac1p}$. Plugging this equality into the second equation, we obtain

\begin{equation} \label{S3.4}
\left\{\begin{array}{rrll}
(-\Delta)^{s} \left( (-\Delta)^{s} u \right)^{\frac1p} &=& u^q & {\rm in} \ \ \Omega\; \\
u &\geq& 0 & {\rm in} \ \ \R^n \\
u &=& 0 & {\rm in} \ \ \R^n\setminus\Omega
\end{array}\right..
\end{equation}

\n The basic idea in trying to solve (\ref{S3.4}) is considering the functional $\Phi : E_p^s \rightarrow \R$ defined by

\begin{equation}\label{func}
\Phi(u)=\frac{p}{p+1}\int\limits_{\Omega}\vert (-\Delta)^{s}u\vert^{\frac{p+1}{p}}dx-\frac{1}{q+1}\int\limits_{\Omega} (u^+)^{q+1}dx\, ,
\end{equation}

\n where $E_p^s$ denotes the reflexive Banach space $W_0^{s, \frac{p+1}{p}}(\Omega) \cap W^{2s, \frac{p+1}{p}}(\Omega)$ as endowed with the norm

\[
\| u \|_{E_p^s} := \left( \int\limits_{\Omega}\vert (-\Delta)^{s}u\vert^{\frac{p+1}{p}}dx \right)^{\frac{p}{p+1}} .
\]

In the case that $E_p^s$ is continuously embedded in $L^{q+1}(\Omega)$, the Gateaux derivative of $\Phi$ at $u \in E_p^s$ in the direction $\varphi \in E_p^s$ is given by

\[
\Phi'(u)\varphi=\int\limits_{\Omega}\left|(-\Delta)^{s}u\right|^{\frac{1}{p}-1} (-\Delta)^{s}u (-\Delta)^{s}\varphi dx - \int\limits_{\Omega}(u^+)^q \varphi dx\, .
\]
This is the case when the couple $(p,q)$ is below the critical hyperbole (\ref{6}).

As we shall see below, weak solutions of (\ref{1}) can be constructed in an appropriate product space from critical points of $\Phi$ in $E_p^s$. Since we do not have a $L^p$ theory available for the fractional Laplace operator, the following result plays an essential role in the proof of the regularity result of the next section.

\begin{propo}\label{prop2}
Let $\Omega$ be a smooth bounded open subset of $\R^{n}$, $n \geq 1$ and $0 < s < 1$. Let $u$ be a critical point of $\Phi$ in $E_p^s$. Assume that the couple $(p,q)$ satisfies (\ref{5}) in the case that $n > 2s$. Then, there exists a function $v \in H_0^{s}(\Omega)$ such that $(u,v)$ is a nonnegative weak solution of the problem (\ref{1}).
\end{propo}
\n {\bf Proof.} Firstly, we claim that $(u^+)^q\in H^{-s}(\Omega)$ for every critical point $u$ of $\Phi$ in $E_p^s$. In fact, note that $L^{(2^{\ast}_s)'}(\Omega)\subset H^{-s}(\Omega)$, where $2^{\ast}_s = \frac{2n}{n-2s}$. So, for $0 < p \leq \frac{2s}{n-2s}$ and $q > 0$, we have $W^{2s, \frac{p+1}{p}}(\Omega) \hookrightarrow L^r(\Omega)$ for all $r \geq 1$ (see \cite{DD, frac}), so that $(u^+)^q\in H^{-s}(\Omega)$. If $p > \frac{2s}{n-2s}$ and $q \leq \frac{n+2s}{n-2s}$, we have $(u^+)^q\in H^{-s}(\Omega)$, because $E_p^s$ is continuously embedded in $L^{q+1}(\Omega)$ (see \cite{DD, frac}) and $\frac{q+1}{q} \geq (2^{\ast}_s)'$. There remains only the cases when $p > \frac{2s}{n-2s}$ and $q>\frac{n+2s}{n-2s}$. By the condition (\ref{5}), we get $p < \frac{n+2s}{n-2s}$. Since $u \in L^{q+1}(\Omega)$ and $q>\frac{n+2s}{n-2s}$, we have $u^+\in H^{-s}(\Omega)$. So, by Proposition 3.4 of \cite{ros131}, the problem

\[
\left\{\begin{array}{rrll}
(-\Delta)^{s} w_1 &=& u^+ & {\rm in} \ \ \Omega\; \\
w_1 &=& 0 & {\rm in} \ \ \R^n\setminus\Omega
\end{array}\right.
\]
admits a unique nonnegative weak solution $w_1 \in H_0^{s}(\Omega)$. For $2s<n\leq 6s$, we have $q+1\geq\frac{n}{2s}$. So, Proposition 1.4 of \cite{ROS} gives $w_1 \in C^{\beta_1}(\overline{\Omega})$ for some $\beta_1\in (0,1)$. But this implies that $u^+\in L^{q(2^{\ast}_s)'}(\Omega)$ (see \cite{S}), so that $(u^+)^q\in H^{-s}(\Omega)$. On the other hand, Sobolev embedding gives us $u^+ \in L^{\frac{n(p+1)}{np-2s(p+1)}}(\Omega)$. Therefore, by Proposition 1.4 of \cite{ROS}, we have $w_1 \in L^{\frac{n(p+1)}{np-4s(p+1)}}(\Omega)$. Since $w_1\in H^{-s}(\Omega)$, again by Proposition 3.4 of \cite{ros131}, the problem

\[
\left\{\begin{array}{rrll}
(-\Delta)^{s} w_2 &=& w_1 & {\rm in} \ \ \Omega\; \\
w_2 &=& 0 & {\rm in} \ \ \R^n\setminus\Omega
\end{array}\right.
\]
admits a unique nonnegative weak solution $w_2 \in H_0^{s}(\Omega)$. If $6s<n\leq 10s$, we have $\frac{n(p+1)}{np-4s(p+1)}\geq\frac{n}{2s}$ and then $w_2\in C^{\beta_2}(\overline{\Omega})$ for some $\beta_2 \in (0,1)$, by Proposition 1.4 of \cite{ROS}. Thus, $w_1$ is a continuous function, and so, $u^+\in L^{q(2^{\ast}_s)'}(\Omega)$ (see \cite{S}). Consequently, $(u^+)^q\in H^{-s}(\Omega)$. Proceeding inductively, we get $(u^+)^q\in H^{-s}(\Omega)$ for any $n>2s$. Using this fact and Proposition 3.4 of \cite{ros131}, it follows that the problem

\[
\left\{\begin{array}{rrll}
(-\Delta)^{s} v &=& (u^+)^q & {\rm in} \ \ \Omega\; \\
v &=& 0 & {\rm in} \ \ \R^n\setminus\Omega
\end{array}\right.
\]
admits a unique nonnegative weak solution $v \in H_0^{s}(\Omega)$. Using the condition (\ref{5}), Proposition 1.4 of \cite{ROS} and Sobolev embedding, we see that $v \in L^{p+1}(\Omega)$. In addition, since $u$ is a critical point of $\Phi$ in $E_p^s$ and $v$ is a weak solution of the above equation in $H_0^{s}(\Omega)$, we have

\begin{equation} \label{7}
\int\limits_{\Omega} \left( \left|(-\Delta)^{s}u\right|^{\frac{1}{p}-1} (-\Delta)^{s}u - v \right) (-\Delta)^{s}\varphi dx = 0
\end{equation}
for all $\varphi \in E_p^s \cap H_0^{s}(\Omega)$ by density (see \cite{DD}). On the other hand, again by Proposition 3.4 of \cite{ros131}, inequality (3.4) in \cite{Francesco} and a consequence of the classical potential theory (see formula (41) of \cite{Stein}), for each $f \in C^\infty_0(\Omega)$ there exists $\varphi \in E_p^s \cap H_0^{s}(\Omega)$ such that

\[
\left\{\begin{array}{rrll}
(-\Delta)^{s} \varphi &=& f & {\rm in} \ \ \Omega\; \\
\varphi &=& 0 & {\rm in} \ \ \R^n\setminus\Omega
\end{array}\right.\, ,
\]
so the equation (\ref{7}) imply that $u$ satisfies in the weak sense

\[
\left\{\begin{array}{rrll}
(-\Delta)^{s} u &=& v^p & {\rm in} \ \ \Omega\; \\
u &=& 0 & {\rm in} \ \ \R^n\setminus\Omega
\end{array}\right.
\]
Lastly, by the weak maximum principle (see \cite{niang}), we conclude that $u$ is nonnegative. In other words, starting from a critical point $u \in E_p^s$ of $\Phi$, we were able to construct a nonnegative weak solution $(u,v)$ of the problem (\ref{1}) in the space $E_p^s \times H_0^{s}(\Omega)$.\ \fim\\

\section{H\"{o}lder regularity}

In this section, we show that weak solutions of (\ref{1}) are H\"{o}lder viscosity solutions by assuming that $(p,q)$ is below the hyperbole (\ref{6}) provided that $n > 2s$.

\begin{propo}\label{prop}
Let $\Omega$ be a smooth bounded open subset of $\R^{n}$, $n \geq 1$ and $0 < s < 1$. Let $(u,v) \in E_p^s \times H_0^{s}(\Omega)$ be a nonnegative weak solution of the problem (\ref{1}). Assume that the couple $(p,q)$ satisfies (\ref{5}) in the case that $n > 2s$, then $(u, v) \in L^{\infty}(\Omega) \times L^{\infty}(\Omega)$ and, moreover, $(u,v) \in C^{\beta}(\R^n) \times C^{\beta}(\R^n)$ for some $\beta \in (0, 1)$.
\end{propo}
\n {\bf Proof.} It suffices to prove the result for $n > 2s$, since the ideas involved in its proof are fairly similar when $n \leq 2s$.

We analyze separately some different cases depending on the values of $p$ and $q$.

For $0 < p < \frac{2s}{n-2s}$ and $q > 0$, we have $W^{2s, \frac{p+1}{p}}(\Omega) \hookrightarrow L^{\infty}(\Omega)$ (see \cite{DD, frac}), so that $u \in L^{\infty}(\Omega)$, and thus, by Proposition 1.4 of \cite{ROS} applied to the second equation of (\ref{1}), one concludes that $v \in L^{\infty}(\Omega)$.

For $\frac{2s}{n-2s} \leq p \leq 1$ and $q>1$, we rewrite the problem (\ref{1}) as follows
\begin{equation}
\left\{
\begin{array}{llll}
(-\Delta)^{s}u = a(x)v^{\frac{p}{2}} & {\rm in} \ \ \Omega\\
(-\Delta)^{s}v = b(x)u & {\rm in} \ \ \Omega\\
u= v=0 & {\rm in} \ \ \R^n\setminus\Omega
\end{array}
\right.,
\end{equation}
where $a(x) = v(x)^{\frac{p}{2}}$ and $b(x) = u(x)^{q - 1}$. Since the couple $(p,q)$ satisfies (\ref{5}), we have $p+1<\frac{2n}{n-2s}$. By Sobolev embedding, $H_0^{s}(\Omega)\rightarrow L^{p+1}(\Omega)$ is bounded, so that $a \in L^{\frac{2(p+1)}{p}}(\Omega)$. Thus, for each fixed $\varepsilon > 0$, we can construct functions $q_{\varepsilon} \in L^{\frac{2(p+1)}{p}}(\Omega)$, $f_{\varepsilon} \in L^{\infty}(\Omega)$ and a constant $K_{\varepsilon} > 0$ such that

\[
a(x) v(x)^{\frac{p}{2}} = q_{\varepsilon}(x)v(x)^{\frac{p}{2}} + f_{\varepsilon}(x)
\]
and

\[
\Vert q_{\varepsilon} \Vert_{L^{\frac{2(p+1)}{p}}(\Omega)} < \varepsilon,\ \ \Vert f_{\varepsilon} \Vert_{L^{\infty}(\Omega)} < K_{\varepsilon}\, .
\]
In fact, consider the set

\[
\Omega_k = \{x \in\Omega: \vert a(x)\vert < k\}\, ,
\]
where $k$ is chosen such that

\[
\int\limits_{\Omega_k^c} \vert a(x)\vert^{\frac{2(p+1)}{p}}dx < \frac{1}{2}\varepsilon^{\frac{2(p+1)}{p}}\, .
\]
This condition is clearly satisfied for $k = k_\varepsilon$ large enough.

We now write

\begin{equation}
q_{\varepsilon}(x)=\left\{
\begin{array}{llll}
\frac{1}{m}a(x) & {\rm for} \ \ x \in \Omega_{k_\varepsilon}\\
a(x) & {\rm for} \ \ x \in \Omega_{k_\varepsilon}^c\\
\end{array}
\right.
\end{equation}
and

\[
f_{\varepsilon}(x) = \left(a(x) - q_{\varepsilon}(x)\right)v(x)^{\frac{p}{2}}\, .
\]
Then,

\begin{eqnarray*}
\int\limits_{\Omega}\vert q_{\varepsilon}(x)\vert^{\frac{2(p+1)}{p}}dx&=&\int\limits_{\Omega_{k_\varepsilon}}\vert q_{\varepsilon}(x)\vert^{\frac{2(p+1)}{p}}dx+\int\limits_{\Omega_{k_\varepsilon}^c}\vert q_{\varepsilon}(x)\vert^{\frac{2(p+1)}{p}}dx\\
&=&\left(\frac{1}{m}\right)^{\frac{2(p+1)}{p}}\int\limits_{\Omega_{k_\varepsilon}}\vert a(x)\vert^{\frac{2(p+1)}{p}}dx+\int\limits_{\Omega_{k_\varepsilon}^c}\vert a(x)\vert^{\frac{2(p+1)}{p}}dx\\
&<&\left(\frac{1}{m}\right)^{\frac{2(p+1)}{p}}\int\limits_{\Omega_{k_\varepsilon}}\vert a(x)\vert^{\frac{2(p+1)}{p}}dx+\frac{1}{2}\varepsilon^{\frac{2(p+1)}{p}}\, .
\end{eqnarray*}
So, for $m = m_\varepsilon > \left(\frac{2^{\frac{p}{2(p+1)}}}{\varepsilon}\right) \Vert a\Vert_{L^{\frac{2(p+1)}{p}}(\Omega)}$, we get

\[
\Vert q_{\varepsilon}\Vert_{L^{\frac{2(p+1)}{p}}(\Omega)} < \varepsilon\, .
\]
Note also that $f_{\varepsilon}(x) = 0$ for all $x \in \Omega_{k_\varepsilon}^c$ and, for this choice of $m$,

\[
f_{\varepsilon}(x) = \left( 1-\frac{1}{m_\varepsilon}\right) a(x)^2 \leq \left( 1-\frac{1}{m_\varepsilon}\right) k_\varepsilon^2
\]
for all $x \in \Omega_{k_\varepsilon}$. Therefore,

\[
\Vert f_{\varepsilon} \Vert_{L^{\infty}(\Omega)} \leq \left( 1-\frac{1}{m_\varepsilon}\right) k_\varepsilon^2 := K_\varepsilon\, .
\]
On the other hand, we have

\[
v(x) = (-\Delta)^{-s}(bu)(x)\, ,
\]
where $b\in L^{\frac{q+1}{q-1}}(\Omega)$. Hence,

\[
u(x) = (-\Delta)^{-s}\left[q_{\varepsilon}(x)((-\Delta)^{-s}(bu)(x))^{\frac{p}{2}}\right] + (-\Delta)^{-s}f_{\varepsilon}(x)\, .
\]

By Proposition 1.4 of \cite{ROS}, the claims $(ii)$ and $(iv)$ below follow readily and, by using H\"{o}lder's inequality, we also get the claims $(i)$ and $(iii)$. Precisely, for fixed $\gamma > 1$, we have:

 \begin{itemize}
    \item[(i)] The map $w \rightarrow b(x)w$ is bounded from $L^\gamma(\Omega)$ to $L^\beta(\Omega)$ for
\[
\frac{1}{\beta} = \frac{q-1}{q+1} + \frac{1}{\gamma};
\]
    \item[(ii)] For any $\theta \geq 1$ in the case that $\beta \geq \frac{n}{2s}$, or for $\theta$ given by
\[
2s=n\left(\frac{1}{\beta}-\frac{2}{p \theta}\right)
\]
in the case that $\beta < \frac{n}{2s}$, there exists a constant $C > 0$, depending on $\beta$ and $\theta$, such that

\[
\Vert ((-\Delta)^sw)^{\frac{p}{2}}\Vert_{L^\theta(\Omega)}\leq C\Vert w\Vert_{L^\beta(\Omega)}^{\frac{p}{2}}
\]
for all $w \in L^\beta(\Omega)$;

    \item[(iii)] The map $ w \rightarrow q_\varepsilon(x)w$ is bounded from $L^\theta(\Omega)$ to $L^\eta(\Omega)$ with norm given by $\Vert q_\varepsilon \Vert_{L^{\frac{2(p+1)}{p}}(\Omega)}$, where $\theta \geq 1$ and $\eta$ satisfies
\[
\frac{1}{\eta}=\frac{p}{2(p+1)}+\frac{1}{\theta};
\]
    \item[(iv)] For any $\delta \geq 1$ in the case that $\eta \geq \frac{n}{2s}$, or for $\delta$ given by

\[
2s = n\left(\frac{1}{\eta} - \frac{1}{\delta}\right)
\]
in the case that $\eta < \frac{n}{2s}$, the map $ w\rightarrow (-\Delta)^{-s}w$ is bounded from $L^\eta(\Omega)$ to $L^\delta(\Omega)$.
    \end{itemize}

Joining $(i)$, $(ii)$, $(iii)$ and $(iv)$ and using that $(p,q)$ satisfies (\ref{5}), one easily checks that $\gamma < \delta$ and, in addition,

\begin{eqnarray*}
\Vert u\Vert_{L^\delta(\Omega)} &\leq & \Vert (-\Delta)^{-s}\left[q_\varepsilon(x)\left((-\Delta)^{-s}(bu)\right)^{\frac{p}{2}}\right]\Vert_{L^\delta(\Omega)}+\Vert (-\Delta)^{-s}f_\varepsilon\Vert_{L^\delta(\Omega)}\\
&\leq & C \left(\Vert q_\varepsilon\Vert_{L^{\frac{2(p+1)}{p}}(\Omega)}\Vert u\Vert_{L^\delta(\Omega)}^{\frac{p}{2}} + \Vert f_\varepsilon\Vert_{L^\delta(\Omega)}\right).
\end{eqnarray*}
Using now the fact that $p \leq 1$, $\Vert q_\varepsilon \Vert_{L^{\frac{2(p+1)}{p}}(\Omega)} < \varepsilon$ and $f_\varepsilon \in L^\infty(\Omega)$, we deduce that $\Vert u \Vert_{L^\delta(\Omega)} \leq C$ for some constant $C > 0$ independent of $u$. Proceeding inductively, we get $u \in L^\delta(\Omega)$ for all $\delta \geq 1$. So, Proposition 1.4 of \cite{ROS} implies that $v \in L^{\infty}(\Omega)$, and thus $u \in L^{\infty}(\Omega)$. Finally, the $C^\beta$ regularity of $u$ and $v$ in $\R^n$ for some $\beta \in (0,1)$ is obtained from each equation by evoking Proposition 1.4 of \cite{ROS}.

The other cases are treated in a similar way by writing $a(x) = v(x)^{p-1}$ if $p > 1$ and $b(x) = u(x)^{\frac{q}{2}}$ if $q \leq 1$ or $b(x) = u(x)^{q-1}$ if $q > 1$.\ \fim\\

\begin{rem} \label{obs}
Thanks to Propositions \ref{prop2} and \ref{prop} and the Silvestre's strong maximum principle (see \cite{S}), we deduce that $(u,v)$ is indeed a positive viscosity solution of (\ref{1}) in $C^0(\R^n) \times C^0(\R^n)$, whenever $u \in E_p^s$ is a nonzero critical point of $\Phi$ and $(p,q)$ is below the critical hyperbole (\ref{6}).
\end{rem}

\section{Rellich variational identity}

In this section, we deduce that positive viscosity solutions of (\ref{1}) satisfy the following integral identity:

\begin{propo}\label{prop1} (Rellich identity)
Let $\Omega$ be a smooth bounded open subset of $\R^n$, $n \geq 1$ and $0 < s < 1$. Then, every positive viscosity solution $(u,v)$ of the problem (\ref{1}) satisfies

\begin{equation} \label{Poho}
\Gamma(1+s)^2\int\limits_{\partial\Omega}\frac{u}{d^s}\frac{v}{d^s} (x\cdot\nu)d\sigma=\left(\frac{n}{q+1} + \frac{n}{p+1} -(n-2s)\right)\int\limits_{\Omega}u^{q+1}dx\, ,
\end{equation}
where $\nu$ denotes the unit outward normal to $\partial\Omega$, $\Gamma$ is the Gamma function, $d(x) = dist(x,\partial\Omega)$ and

\[
\frac{u}{d^s}(x) := \lim_{\varepsilon \rightarrow 0^+} \frac{u(x - \varepsilon \nu)}{d^s(x - \varepsilon \nu)} > 0
\]
for all $x \in \partial \Omega$.
\end{propo}

It deserves mention that $u/d^s, v/d^s \in C^\alpha(\overline{\Omega})$ for some $\alpha \in (0,1)$ and $u/d^s,v/d^s>0$ in $\overline{\Omega}$ (see Theorem 1.2 in \cite{ROS1} or Proposition 2.7 in \cite{CR}). So, the left-hand side of the identity (\ref{Poho}) is well defined.

\n {\bf Proof.} Let $(u,v)$ be a viscosity solution of (\ref{1}). Then,

\begin{equation}
\left\{
\begin{array}{rrll}
(-\Delta)^{s} (u+v) &=& v^p+u^q & {\rm in} \ \ \Omega\\
u+v &=& 0  & {\rm in} \ \ \R^n\setminus\Omega
\end{array}
\right.
\end{equation}
and
\begin{equation}
\left\{
\begin{array}{rrll}
(-\Delta)^{s} (u-v) &=& v^p-u^q & {\rm in} \ \ \Omega\\
u-v &=& 0  & {\rm in} \ \ \R^n\setminus\Omega
\end{array}
\right.
\end{equation}

Applying the Pohozaev variational identity for semilinear problems involving the operator $(-\Delta)^{s}$ (Theorem 1.1 of \cite{ROS2}), we get

\begin{eqnarray*}
-\int\limits_{\Omega}(x\cdot\nabla u+v)((-\Delta)^{s}u+(-\Delta)^{s}v)dx&=&-\frac{2s-n}{2}\int\limits_{\Omega}(u+v)(v^p+u^q)dx\\
& &+\frac{1}{2}\Gamma(1+s)^2\int\limits_{\partial\Omega}\left(\frac{u+v}{d^s}\right)^2(x\cdot\nu)d\sigma
\end{eqnarray*}
and
\begin{eqnarray*}
-\int\limits_{\Omega}(x\cdot\nabla u+v)((-\Delta)^{s}u-(-\Delta)^{s}v)dx&=&-\frac{2s-n}{2}\int\limits_{\Omega}(u-v)(v^p-u^q)dx\\
& &+\frac{1}{2}\Gamma(1+s)^2\int\limits_{\partial\Omega}\left(\frac{u-v}{d^s}\right)^2(x\cdot\nu)d\sigma\, .
\end{eqnarray*}
Now subtracting both identities, we obtain

\begin{eqnarray}\nonumber
2\int\limits_{\Omega}[(x\cdot\nabla u)(-\Delta)^{s}v+(x\cdot\nabla v)(-\Delta)^{s}u]dx&=&(2s-n)\int\limits_{\Omega}[u(-\Delta)^{s}v+v(-\Delta)^{s}u]dx\\\label{5.9}
& &-2\Gamma(1+s)^2\int\limits_{\partial\Omega}\frac{u}{d^s}\frac{v}{d^s} (x\cdot\nu)d\sigma\, .
\end{eqnarray}
Because $v=0$ in $\R^n\setminus\Omega$, we have
\[
\int\limits_{\Omega}(x\cdot\nabla u)(-\Delta)^{s}vdx=\int\limits_{\Omega}(x\cdot\nabla u)u^qdx=\frac{1}{q+1}\int\limits_{\Omega}(x\cdot\nabla u^{q+1})dx=-\frac{n}{q+1}\int\limits_{\Omega}u^{q+1}dx\, .
\]
In a similar way,
\[
\int\limits_{\Omega}(x\cdot\nabla v)(-\Delta)^{s}udx=-\frac{n}{p+1}\int\limits_{\Omega}v^{p+1}dx\, .
\]
Plugging these two identities into (\ref{5.9}), we derive
\[
2\Gamma(1+s)^2\int\limits_{\partial\Omega}\frac{u}{d^s}\frac{v}{d^s} (x\cdot\nu)d\sigma=\left(2s-n+\frac{2n}{q+1}\right)\int\limits_{\Omega}u^{q+1}dx+\left(2s-n+\frac{2n}{p+1}\right)\int_{\Omega}v^{p+1}dx\, .
\]
Since every viscosity solution of (\ref{1}) is also a bounded weak solution, one has
\[
\int\limits_{\Omega}v^{p+1} dx=\int\limits_{\Omega}v(-\Delta)^{s}u dx=\int\limits_{\R^n}(-\Delta)^{s/2}u(-\Delta)^{s/2}v dx=\int\limits_{\Omega}u(-\Delta)^{s}v dx=\int\limits_{\Omega}u^{q+1} dx\, .
\]
Thus, the desired conclusion follows directly from this equality.\; \fim\\

\section{Proof of Theorem \ref{teo1}}

We organize the proof of Theorem \ref{teo1} into two parts. We start by proving the existence of a positive viscosity solution. By Remark \ref{obs}, it suffices to show the existence of a nonzero critical point $u \in E^s_p$ of the functional $\Phi$.

\subsection{The existence part}
We apply the direct method to the functional $\Phi$ on $E^s_p$.

In order to show the coercivity of $\Phi$, note that $q + 1 < \frac{p+1}{p}$ because $pq < 1$, so that the embedding $E^s_p \hookrightarrow L^{q+1}(\Omega)$ is continuous. So, there exist constants $C_1, C_2 > 0$ such that
\begin{eqnarray*}
\Phi(u) &=& \frac{p}{p+1}\int\limits_{\Omega}\vert (-\Delta)^{s}u\vert^{\frac{p+1}{p}}dx - \frac{1}{q+1}\int\limits_{\Omega}\vert u\vert ^{q+1}dx\\
&\geq& \frac{p C_1}{p+1}\Vert u\Vert_{E^s_p}^{\frac{p+1}{p}}-\frac{C_2}{q+1}\Vert u\Vert_{E^s_p}^{q+1}\\
&=& \Vert u\Vert^{\frac{p+1}{p}}_{E^s_p}\left(\frac{p C_1}{p+1}-\frac{C_2}{(q+1)\Vert u\Vert^{\frac{p+1}{p}-(q+1)}_{E^s_p}}\right)
\end{eqnarray*}
for all $u\in E^s_p$. For the existence of $C_1$ see \cite{DongKim, Francesco, Zlemer, Stein}. Therefore, $\Phi$ is lower bounded and coercive, that is, $\Phi(u) \rightarrow +\infty$ as $\Vert u\Vert_{E^s_p} \rightarrow +\infty$.

Let $(u_k) \subset E^s_p$ be a minimizing sequence of $\Phi$. It is clear that $(u_k)$ is bounded in $E^s_p$, since $\Phi$ is coercive. So, module a subsequence, we have $u_k \rightharpoonup u_0$ in $E^s_p$. Since $E^s_p$ is compactly embedded in $L^{q+1}(\Omega)$ (see \cite{DD, frac}), we have $u_k \rightarrow u_0$ in $L^{q+1}(\Omega)$. Here, we again use the fact that $q + 1 < \frac{p+1}{p}$. Thus,

\begin{eqnarray*}
\lim_{n\rightarrow\infty}\inf \Phi(u_k) &=& \lim_{k \rightarrow \infty} \inf \frac{p}{p+1}\Vert (-\Delta)^{s} u_k \Vert^{\frac{p+1}{p}}_{L^{\frac{p+1}{p}}(\Omega)}-\frac{1}{q+1}\Vert u_{0}\Vert_{L^{q+1}(\Omega)}^{q+1}\\
&\geq &\frac{p}{p+1}\Vert (-\Delta)^{s}u_{0}\Vert^{\frac{p+1}{p}}_{L^{\frac{p+1}{p}}(\Omega)}-\frac{1}{q+1}\Vert u_{0}\Vert_{L^{q+1}(\Omega)}^{q+1} = \Phi(u_{0})\, ,
\end{eqnarray*}
so that $u_0$ minimizers $\Phi$ on $E^s_p$. We just need to guarantee that $u_0$ is nonzero. But, this fact is clearly true since $\Phi(\varepsilon u_1) < 0$ for any nonzero nonnegative function $u_1 \in E^s_p$ and $\varepsilon > 0$ small enough, that is,

\[
\Phi(\varepsilon u_1) = \frac{p\varepsilon^{\frac{p+1}{p}}}{p+1}\int\limits_{\Omega}\vert (-\Delta)^{s} u_1 \vert^{\frac{p+1}{p}}dx - \frac{\varepsilon^{q+1}}{q+1}\int\limits_{\Omega}\vert u_1 \vert ^{q+1}dx < 0
\]
for $\varepsilon>0$ small enough. This ends the proof of existence.\; \fim\\

\subsection{The uniqueness part}
The main tools in the proof of uniqueness are the Silvestre's strong maximum principle, a $C^\alpha$ regularity result up to the boundary and a Hopf's lemma adapted to fractional operators.

Let $(u_{1},v_{1}),(u_{2},v_{2}) \in C^0(\R^n) \times C^0(\R^n)$ be two positive viscosity solutions of (\ref{1}). Define
\[
S=\{s\in(0,1] : u_{1} - tu_{2},\ v_{1} - tv_{2}\geq 0\ \text{ in } \overline{\Omega}\ \text{ for all }\ t \in [0,s]\}\, .
\]
By Theorem 1.2 of \cite{ROS1} (or Proposition 2.7 of \cite{CR}), we have $u_i/d^s, v_i/d^s \in C^\alpha(\overline{\Omega})$ and both quotients are positive on $\overline{\Omega}$. So, $(u_{1}-tu_{2})/d^s, (v_{1}-tv_{2})/d^s > 0$ on $\partial \Omega$ for $t > 0$ small enough and thus the set $S$ is no empty.

Let $s_{\ast}=\sup S$ and assume that $s_{\ast} < 1$.

Clearly,

\begin{equation}\label{eq1}
u_{1}-s_{\ast}u_{2},\ v_{1}-s_{\ast}v_{2}\geq 0\ \text{ in }\ \overline{\Omega}\, .
\end{equation}
By (\ref{eq1}) and the integral representation in terms of the Green function $G_{\Omega}$ of $(-\Delta)^{s}$ (see \cite{Green, ros176}), we have

\begin{eqnarray*}
u_{1}(x)&=&\int\limits_{\Omega}G_{\Omega}(x,y)v_{1}^{p}(y)dy\geq\int\limits_{\Omega}G_{\Omega}(x,y)s_{\ast}^p v_{2}^{p}(y)dy\\
&=& s_{\ast}^{p}\int\limits_{\Omega}G_{\Omega}(x,y)v_{2}^{p}(y)dy=s_{\ast}^{p}u_{2}(x)
\end{eqnarray*}
for all $x \in \overline{\Omega}$. In a similar way, one gets $v_1 \geq s_{\ast}^{q} v_2$ in $\overline{\Omega}$.

\n Using the assumption $pq < 1$ and the fact that $s_{\ast} < 1$, we derive

\begin{equation}
\left\{\begin{array}{llll}
(-\Delta)^{s}(u_{1} - s_{\ast}u_{2}) = v_{1}^{p}-s_{\ast}v_{2}^{p} \geq (s_{\ast}^{pq}-s_{\ast})v_{2}^{p} > 0  \\
(-\Delta)^{s}(v_{1} - s_{\ast}v_{2}) = u_{1}^{q}-s_{\ast}u_{2}^{q} \geq (s_{\ast}^{pq}-s_{\ast})u_{2}^{q} > 0 \ \
\end{array}\right.\ {\rm in}\ \Omega
\end{equation}
So, by the Silvestre's strong maximum principle (see \cite{S}), one has $u_{1}-s_{\ast}u_{2}, v_{1}-s_{\ast}v_{2} > 0$ in $\Omega$. Again, arguing as above, we easily deduce that $(u_{1}-s_{\ast}u_{2})/d^s,(v_{1}-s_{\ast}v_{2})/d^s>0$ on $\partial \Omega$, so that $u_{1} - (s_{\ast}+\varepsilon)u_{2}, v_{1} - (s_{\ast}+\varepsilon)v_{2} > 0$ in $\Omega$ for $\varepsilon > 0$ small enough, contradicting the definition of $s_{\ast}$. Therefore, $s_{\ast} \geq 1$ and, by (\ref{eq1}), $u_{1} - u_{2}, v_{1} - v_{2} \geq 0$ in $\overline{\Omega}$. A similar reasoning also produces $u_{2} - u_{1}, v_{2} - v_{1} \geq 0$ in $\overline{\Omega}$. This ends the proof of uniqueness.\; \fim\\

\section{Proof of Theorem \ref{teo2}}

By Remark \ref{obs}, it suffices to show the existence of a nonzero critical point $u \in E^s_p$ of the functional $\Phi$. Assume $p, q > 0$, $pq > 1$ and the assumption (\ref{5}). The proof consists in applying the classical mountain pass theorem of Ambrosetti and Rabinowitz in our variational setting. By well-known embedding theorems (see \cite{DD, frac}), (\ref{5}) implies that $E^s_p$ is compactly embedded in $L^{q+1}(\Omega)$. We first assert that $\Phi$ has a local minimum in the origin. Consider the set $\Gamma :=\left\{u\in E^s_p : \Vert u\Vert_{E^s_p}=\rho\right\}$. Then, on $\Gamma$, we have

\begin{eqnarray*}
\Phi(u) &=& \frac{p}{p+1}\int\limits_{\Omega}\vert (-\Delta)^{s}u\vert^{\frac{p+1}{p}}dx-\frac{1}{q+1}\int\limits_{\Omega}\vert u\vert^{q+1}dx\\
&\geq & \frac{p C_1}{p+1}\Vert u\Vert_{E^s_p}^{\frac{p+1}{p}} - \frac{C_2}{q+1}\Vert u\Vert_{E^s_p}^{q+1} =\rho^{\frac{p+1}{p}}\left( \frac{p C_1}{p+1}-\frac{C_2}{q+1}\rho^{q+1-\frac{p+1}{p}}\right) \\
&>& 0 = \Phi(0)
\end{eqnarray*}
for fixed $\rho>0$ small enough, so that the origin $u_0 = 0$ is a local minimum point. For the existence of $C_1$ see \cite{DongKim, Francesco, Zlemer, Stein}. In particular,  $\inf_{\Gamma} \Phi > 0 = \Phi(u_0)$.

Note that $\Gamma$ is a closed subset of $E^s_p$ and decomposes $E^s_p$ into two connected components, namely $\left\{u\in E^s_p : \Vert u\Vert_{E^s_p} < \rho \right\}$ and $\left\{u\in E^s_p : \Vert u\Vert_{E^s_p} > \rho\right\}$.

Let $u_1 = t \overline{u}$, where $t > 0$ and $\overline{u} \in E^s_p$ is a nonzero nonnegative function. Since $pq > 1$, we can choose $t$ sufficiently large so that
\[
\Phi(u_1)=\frac{p t^{\frac{p+1}{p}}}{p+1} \int\limits_{\Omega} \vert (-\Delta)^{s} \overline{u} \vert^{\frac{p+1}{p}}dx - \frac{t^{q+1}}{q+1} \int\limits_{\Omega} (\overline{u}^+)^{q+1}dx < 0\, .
\]
It is clear that $u_1 \in \left\{u\in E^s_p: \Vert u\Vert_{E^s_p} > \rho\right\}$. Moreover, $\inf_{\Gamma} \Phi > \max\{\Phi(u_0), \Phi(u_1)\}$, so that the mountain pass geometry is satisfied.

Finally, we show that $\Phi$ fulfills the Palais-Smale condition (PS). Let $(u_k) \subset E^s_p$ be a (PS)-sequence, that is,
\[
\vert \Phi(u_k)\vert \leq C_0
\]
and

\[
\vert \Phi'(u_k) \varphi\vert \leq \varepsilon_k \Vert \varphi \Vert_{E^s_p}
\]
for all $\varphi \in E^s_p$, where $\varepsilon_k \rightarrow 0$ as $k \rightarrow +\infty$.

From these two inequalities and the assumption $pq > 1$, we deduce that

\begin{eqnarray*}
C_0 + \varepsilon_k \Vert u_k \Vert_{E^s_p} &\geq & \vert (q+1) \Phi(u_k) - \Phi'(u_k) u_k\vert \\
&\geq &\left(\frac{p(q+1)}{p+1}-1\right) \int\limits_{\Omega}\vert (-\Delta)^{s}u_k\vert^{\frac{p+1}{p}}dx\\
&\geq & C \Vert u_k\Vert_{E^s_p}^{\frac{p+1}{p}}
\end{eqnarray*}
and thus $(u_k)$ is bounded in $E^s_p$. Thanks to the compactness of the embedding $E^s_p \hookrightarrow L^{q+1}(\Omega)$, one easily checks that $(u_k)$ converges strongly in $E^s_p$. So, by the mountain pass theorem, we obtain a nonzero critical point $u \in E^s_p$. This ends the proof.\\

\section{Proof of Theorem \ref{teo3}}

It suffices to assume that $\Omega$ is star-shaped with respect to the origin, that is, $(x\cdot\nu) > 0$ for any $x\in\partial\Omega$, where $\nu$ is the unit outward normal to $\partial\Omega$ at $x$.

Arguing by contradiction, assume the problem (\ref{1}) admits a positive viscosity solution $(u,v)$. Then, by Theorem 1.2 of \cite{ROS1}, we have

\[
2\Gamma(1+s)^2\int\limits_{\partial\Omega}\frac{u}{d^s}\frac{v}{d^s} (x\cdot\nu)d\sigma > 0\, .
\]

On the other hand, the assumption (\ref{4}) is equivalent to $\frac{n}{q+1} + \frac{n}{p+1} -(n-2s) \leq 0$ and thus the right-hand side of the Rellich identity (\ref{Poho}) is non-positive, providing the wished contradiction. Hence, the problem (\ref{1}) admits no positive viscosity solution and we end the proof.\; \fim\\

\n {\bf Acknowledgments:} The first author was supported by CAPES and the second one was supported by CAPES (BEX 6961/14-2), CNPq (PQ 306406/2013-6) and Fapemig (PPM 00223-13). \\

\end{document}